\documentstyle{article}
\title{\centerline{SYMMETRY GROUPS AND LAGRANGIANS}
\centerline{ASSOCIATED TO \c{T}I\c{T}EICA SURFACES}}
\author{N.B\^IL\u{A}}
\date{}
\begin{document}
\maketitle

\begin{abstract}

\c{T}i\c{t}eica proved [23] that the surfaces for which the ratio
${K\over {d^4}}$ is constant (where $K$ is the Gaussian curvature and $d$ is
the distance from the origin to the tangent plane at an arbitrary point) are
invariants under the group of centroaffine transformations. In this paper
one applies the symmetry groups theory ([20], [21]) to study the PDEs 
which arise in \c{T}i\c{t}eica surfaces theory:
the PDEs system (5) (equivalent to (15)+(16)) with the particular cases
(7) and (9), the PDEs system (11) and the Liouville-\c{T}i\c{t}eica
PDE (8), respectively the \c{T}i\c{t}eica PDE (10) (equivalent to (8'),
respectively (10')). In the case of the PDEs systems (15)+(16) and (11),
the center of our attention is to find the symmetry subgroups $G_1$ of
the full symmetry group $G$, respectively $\bar G_1$
of the full symmetry group $\bar G$,
which act on the space of dependent variables, and also
the symmetry subgroups $G_2$ of $G$ and $\bar G_2$ of $\bar G$, which
act on the space of independent variables (Theorems 4,5,6 and 7).
One proves (Theorem 4) that the subgroup $G_1$ is the unimodular
subgroup of the group of the centroaffine transformations.
One gives a new solution (21) for the \c{T}i\c{t}eica PDE
for which it is a ruled \c{T}i\c{t}eica surface
associated (Proposition 1).
One finds the symmetry groups of Liouville-\c{T}i\c{t}eica PDE
and \c{T}i\c{t}eica PDE (Theorems 9 and 10)
and one proves that these are Euler-Lagrange equations 
with the Lagrangians (32) and (33) (Theorem 12).
One gets the variational symmetry groups of the associated functionals
(38) and (39) (Theorems 16 and 17) and conservations laws (Proposition 3).
We make the remark that the \c{T}i\c{t}eica simple surfaces $z=f(x,y)$
was studied also in [28] and we proved that the \c{T}i\c{t}eica PDE
is an Euler-Lagrange equation.
All these results shows that \c{T}i\c{t}eica surfaces theory
is strongly related to variational problems and hence it is a subject
of global differential geometry.

\end{abstract}


\section{Introduction}

\c{T}i\c{t}eica-the founder of the centroaffine geometry- introduced in 1907
a new class of surfaces, called {\it the surfaces S}, with the property that
${K\over {d^4}}$=constant, where $K$ is the Gaussian curvature and $d$ is the
distance from the origin to the tangent plane at an arbitrary point
[23]. These were called {\it \c{T}i\c{t}eica surfaces}
by Gheorghiu, or {\it affine spheres} by Blaschke and {\it projectives spheres}
by Wilczynski.
The most simple \c{T}i\c{t}eica surfaces are the spheres and the quadrics.
The extension of this class to hypersurfaces was considered by
\c{T}i\c{t}eica itself [24],[25]. Also, Mayer [19], Gheorghiu [12],
Dobrescu [10] and Vr\u{a}nceanu [29], [30] studied the properties of these
hypersurfaces. Gheorghiu maked a remark on the hypersurfaces \c{T}i\c{t}eica:
these can be considered as the affine spaces $A_{n-1}$, embedded in a affine
Euclidean space $E_n$.
Using this result, he introduced a new class of affine space $A^0_n$ and new
examples of these were considered by Udri\c{s}te [26].
We start to make a short presentation of the \c{T}i\c{t}eica surfaces.

Let $D\subset ${\bf R}$^2$ be an open set and let
$$
\Sigma :\;\;{\bf r} (u,v)=x(u,v){\bf i}+y(u,v){\bf j}+z(u,v){\bf k},\;\;
(u,v)\in D,
$$
be a surface in {\bf R}$^3$, different from a cone with the vertex at the
origin of the system of coordinates.
Thus, the position vector {\bf r} of an arbitrary point of the surface satisfies the
condition
$$
({\bf r},{\bf r}_u,{\bf r}_v)\neq 0, \leqno(1)
$$
and this can be considered the solution of the second order PDEs system,
$$
\left\{
\begin{array}{ccl}
{\bf r}_{uu} & = & a{\bf r}_u+b{\bf r}_v+c{\bf r} \\
{\bf r}_{uv} & = & a'{\bf r}_u+b'{\bf r}_v+c'{\bf r} \\
{\bf r}_{vv} & = & a''{\bf r}_u+b''{\bf r}_v+c''{\bf r}, \\
\end{array} \right. \leqno(2)
$$
which is completely integrable, i.e.,
$$
({\bf r}_{uu})_v = ({\bf r}_{uv})_u,\;\;\;
({\bf r}_{uv})_v =  ({\bf r}_{vv})_u, \leqno(3)
$$
where $a,a',a'',...$ are nine functions of $u$ and $v$.
The above system defines a surface, leaving a centroaffinity aside.
The coeficients $a,a',a'',...$ are called {\it the centroaffine invariants}.
If the surface $\Sigma $ is related to the asymptotic lines (if a surface
is not developable, then the two families of the asymptotic lines are
distinct), then $c=c''=0$, and thus it is defined by the following completely
integrable second order PDEs system
$$
\left\{
\begin{array}{ccl}
{\bf r}_{uu} & = & a{\bf r}_u+b{\bf r}_v \\
{\bf r}_{uv} & = & a'{\bf r}_u+b'{\bf r}_v+c'{\bf r} \\
{\bf r}_{vv} & = & a''{\bf r}_u+b''{\bf r}_v. \\
\end{array} \right. \leqno(4)
$$

{\bf Theorem 1 (\c{T}i\c{t}eica)}.
{\it Let $\Sigma $ be a surface related to the asymptotic lines. The ratio
$I={K\over {d^4}}$ is a constant if and only if} $a'=b'=0$.

Thus, the surfaces \c{T}i\c{t}eica are defined by the PDEs system
$$
\left\{
\begin{array}{ccl}
{\bf r}_{uu} & = & a{\bf r}_u+b{\bf r}_v \\
{\bf r}_{uv} & = & h{\bf r} \\
{\bf r}_{vv} & = & a''{\bf r}_u+b''{\bf r}_v, \\
\end{array} \right. \leqno(5)
$$
where we denote $c'=h$, and for which the integrability conditions (3) turn in
$$
ah=h_u,\;\;a_v=ba''+h,\;\;b_v+bb''=0, \leqno(6)
$$
$$
h_v=b''h,\;\;a''_u+aa''=0,\;\;h=b''_u+a''b.
$$

{\bf Remark}. In the particular cases $b=0$ or $a''=0$, $\Sigma $ is a simply
ruled surface. Thus, for $b=0,\;a''\neq 0$: the coordinates curves $v=v_0$ are
straight lines, and for $b\neq 0,\;a''=0$: the coordinates curves $u=u_0$ are
straight lines. If $b=a''=0$, then $\Sigma $ is a double ruled surface (a
quadric surface).

{\it The ruled \c{T}i\c{t}eica surfaces} are given by the PDEs system
$$
\left\{
\begin{array}{ccl}
{\bf r}_{uu} & = & {h_u\over {h}}{\bf r}_u+{\varphi (u)\over {h}}{\bf r}_v \\
{\bf r}_{uv} & = & h{\bf r} \\
{\bf r}_{vv} & = & {h_v\over {h}}{\bf r}_v, \\
\end{array} \right. \leqno(7)
$$
where $h$ is a solution of the Liouville-\c{T}i\c{t}eica PDE
$$
(\ln h)_{uv}=h. \leqno(8)
$$

{\it The \c{T}i\c{t}eica surfaces which are not ruled surfaces}, are given by
the PDEs system
$$
\left\{
\begin{array}{ccl}
{\bf r}_{uu} & = & {h_u\over {h}}{\bf r}_u+{1\over {h}}{\bf r}_v \\
{\bf r}_{uv} & = & h{\bf r} \\
{\bf r}_{vv} & = & {1\over {h}}{\bf r}_u+{h_v\over {h}}{\bf r}_v, \\
\end{array} \right. \leqno(9)
$$
where $h$ is a solution of the \c{T}i\c{t}eica PDE
$$
(\ln h)_{uv}=h-{1\over {h^2}}. \leqno(10)
$$

The PDEs system (5) can be identified with the completely integrable scalar
PDEs system
$$
\left\{
\begin{array}{ccl}
\theta _{uu} & = & a\theta _u+b\theta _v \\
\theta _{uv} & = & h\theta  \\
\theta _{vv} & = & a''\theta _u+b''\theta_v, \\
\end{array} \right.   \leqno(11)
$$
with the condition that three independent solutions of (11)+(6):
$x=x(u,v),y=y(u,v),z=z(u,v),$ define a \c{T}i\c{t}eica surface.
It is known that every linear combination of $x,y,z$ is a solution of the
system (11) also. Thus, a surface $\Sigma $ is determined leaving a
centroaffinity aside.

On the other hand, it is known that Sophus Lie is the founder of theory of
symmetry groups.
A modern presentation using the jets theory is introduced by Olver
in his book [20].
A local group of transformations on the
space of the independent and dependent variables of a studied PDEs system
which transforms the solutions of the system into its solutions,
is called {\it symmetry group} or {\it strong symmetry group}
of the system.
The symmetry groups theory is very applied to study the ODEs, PDEs
systems which appear in Geometry, Mechanics and Physics
[2]-[6],[8],[9],[14],[17],[18],[20],[21],[27],[28],[31].
There are many computational programs for finding
the defining system of infinitesimal symmetries, but
in the ours cases of the PDEs systems (5) and (11) we cannot apply these.
We make the remark that a other point of view in the study of
the Liouville-\c{T}i\c{t}eica PDE
and \c{T}i\c{t}eica PDE is contained the papers of
Bobenko [7] and Wolf [31].

In this paper we shall apply this theory for finding
infinitesimal symmetries of the PDEs systems which arise
in \c{T}i\c{t}eica surfaces theory
and we shall give a new point of view of \c{T}i\c{t}eica theory
with the connection of the known results.
We shall adopt the notation of the book of Olver [20].

\section{Symmetry Group of PDEs System}

Let consider the PDEs system
$$
\Delta _\nu (x,{u}^{(n)})=0,\;\;\nu =1,...,l, \leqno(12)
$$
with $x=(x^1,...,x^p),\;u=(u^1,...,u^q)$ and
$\Delta (x,u^{(n)})=(\Delta _1(x,u^{(n)}),...,\Delta _l(x,u^{(n)})),$
is a differentiable function.
We note $u^{(n)}$ all the partial derivatives of the function
$u$ to 0 up $n$. Any function $u=h(x)$,
$h:D\subset {\bf R}^p\to U\subset {\bf R}^q,\;\;h=(h^1,...,h^q),$
induces the function $u^{(n)}=pr^{(n)}h$ called {\it the
n-th prolongation of h}, which is defined by $u^\alpha _J=
\partial _Jh^\alpha $, $pr^{(n)}h:D\to U^{(n)}$,
and for each $x\in D$, $pr^{(n)}h$ is a vector whose $qp^{(n)}=C^{n}_{p+n}$
entries represent the values of $h$ and all its derivatives up to order
$n$ at the point $x$.

The space $D\times U^{(n)}$, whose coordinates represent the
independent variables, the dependent variables and the derivatives of the
dependent variables up to order $n$, is called {\it the n-th order
jet space} of the underlying space $D\times U$.
Thus $\Delta $ is a map from the jet space $D\times U^{(n)}$ to ${\bf R}^l$.
The PDEs system (12) determine the subvariety
$$
{\cal S}=\{(x,u^{(n)})\;\vert \Delta (x,u^{(n)})=0\}
$$
of the total jet space $D\times U^{(n)}$. One identifies the system of PDEs
(12) with its corresponding subvariety ${\cal S}$.

Let $M\subset D\times U$ be an open set.
{\it A symmety group of the PDEs system (12)}
is a local group of transformations $G$ acting on $M$
with the property that whenever $u=f(x)$ is a solution of (12)
and whenever $g\cdot f$ is defined for $g\in G$, then $u=g\cdot f(x)$
is also a solution of the system.
The system (12) is called {\it invariant with respect to $G$}.

Let us consider $X$ a vector field on $M$ 
with corresponding (local) 1-parameter group $\exp (\varepsilon X)$
which is the infinitesimal generator of the symmetry group of the
PDEs system (12).
The infinitesimal generator
of the corresponding prolonged 1-parameter group
$pr^{(n)}[\exp (\varepsilon X)]$:
$$
pr^{(n)}X\vert _{(x,u^{(n)})}={d\over {d\varepsilon }}
pr^{(n)}[\exp(\varepsilon X)](x,u^{(n)})\vert _{\varepsilon =0}
$$
for any $(x,u^{(n)})\in M^{(n)}$,
is a vector field on the $n$-jet space $M^{(n)}$ called
{\it the n-th prolongation of X} and denoted by $pr^{(n)}X$.

The PDEs system (12) is called to be of {\it maximal rank} if the Jacobi
matrix
$$
J_\Delta (x,u^{(n)})=\left({\partial \Delta _\nu
\over {\partial x^i}},{\partial \Delta _{\nu }\over {\partial u^\alpha _J}}
\right)
$$
of $\Delta $, with respect to all the variables $(x,u^{(n)})$,is of rank $l$
whenever $\Delta (x,u^{(n)})=0$.

{\bf Theorem 2}.
{\it Let
$$
X=\sum _{i=1}^p\zeta ^i(x,u){\partial \over {\partial x^i}}+
\sum _{\alpha =1}^q\phi _{\alpha }(x,u){\partial \over {\partial u^\alpha }}
$$
be a vector field on open set $M\subset D\times U$. The $n$-th prolongation of
$X$ is the vector field
$$
pr^{(n)}X=X+\sum _{\alpha =1}^q\sum _J \phi _{\alpha }^J(x,
u^{(n)}){\partial \over {\partial u^{\alpha }_J}}, \leqno(13)
$$
defined on the corresponding jet space $M^{(n)}\subset D\times U^{(n)}$,
the second summation being over all multi-indices
$J=(j_1,...j_k)$ whith $1\leq j_k\leq p$, $1\leq k\leq n$.
The coefficient functions $\phi ^J_{\alpha }$ of $pr^{(n)}X$ are given by the
following formula
$$
\phi ^J_{\alpha }(x,u^{(n)})=D_J\left(\phi _\alpha -
\sum _{i=1}^p\zeta ^iu^\alpha _i\right)+\sum _{i=1}^p\zeta ^iu_{J,i}^\alpha ,
$$
where $u^\alpha _i={\partial u^{\alpha } \over {\partial x^i}},\;\;
u^\alpha _{J,i}={\partial u^{\alpha }_J \over {\partial x^i}}.$}

{\bf Theorem 3 (Criterion of infinitesimal invariance)}.
{\it Let us consider the PDEs system (12) of maximal rank defined over $M\subset
D\times U$. If $G$ is a local group of transformations acting on $M$ and
$$
pr^{(n)}X[\Delta _{\nu }(x,u^n)]=0,\;\;\nu =1,...,l,\leqno(14)
$$
whenever $\Delta  _\nu (x,u^{(n)})=0$, for every infinitesimal generator $X$ of $G$,
then $G$ is a symmetry group of the PDEs system (12).}

{\bf The algorithm for finding the symmetry group G of the
PDEs system (12)}:
One considers the vector field $X$ on $M$ and one writes the infinitesimal
invariance condition (14);
one eliminates any dependence between partial derivatives of the functions
$u^\alpha $, determined by the PDEs system (12);
one writes the condition (14) like polynomials in the partial derivatives
of $u^{\alpha }$;
one equates with zero the coefficients of partial derivatives of $u^{\alpha }$
in (14);
it follows a PDEs system with respect to the unknown functions
$\zeta ^i,\; \phi _{\alpha }$ and this system defines the symmetry group $G$
of the studied PDEs system.

\section{Symmetry Groups
Associated to PDEs Systems of \c{T}i\c{t}eica Surfaces}

{\bf 3.1}.
In the first part of this section we shall study the symmetries of the PDEs
system (5), which can be considered in the equivalent form
$$
\left\{
\begin{array}{ccl}
x_{uu} & = & ax_u+bx_v \\
x_{uv} & = & hx \\
x_{vv} & = & a''x_u+b''x_v \\
y_{uu} & = & ay_u+by_v   \\
y_{uv} & = & hy   \\
y_{vv} & = & a''y_u+b''y_v \\
z_{uu} & = & az_u+bz_v   \\
z_{uv} & = & hz     \\
z_{vv} & = & a''z_u+b''z_v, \\
\end{array} \right. \leqno(15)
$$
with the conditions (1) and (6). The condition (1) can be written as
$$
(y_uz_v-z_uy_v)x-(x_uz_v-x_vz_u)y+(x_uy_v-x_vy_u)z=f , \leqno(16)
$$
where $f=f(u,v)$ is a nonzero function.
We consider the case of the real asymptotic lines.
Let $D\times U^{(2)}$ be the second order jet space associated to the
PDEs system (15)+(16), whose coordinates represent the independent variables $u,v$,
the dependent variables $x,y,z$ and the derivatives of the dependent variables
till the order two.
Denote $x^1=u,\;x^2=v$, $u^1=x,\;u^2=y$ and $u^3=z$ (in the above section).
Let $M\subset D\times U$ be an open set and let
$$
X=\zeta {\partial \over {\partial u}}+\eta {\partial \over {\partial v}}
+\phi {\partial \over {\partial x}}+\lambda {\partial \over {\partial y}}+
\psi {\partial \over {\partial z}}
$$
be the infinitesimal generator of the symmetry group
$G$ of the PDEs system (15)+(16),
where $\zeta ,\eta ,\phi ,\lambda ,\psi $ are functions of $u,v,x,y$ and $z$.

We shall study if there is a subgroup $G_1$ of the
symmetry group $G$, which acts
on the space of the dependent variables $x,y,z$ of the given system.
Let suppose that it is, and let $Y$ be its infinitesimal generator.
In this case, we must have
$\zeta = 0,\;\eta = 0,\;\phi =\phi (x,y,z),\;\lambda =\lambda (x,y,z),\;
\psi =\psi (x,y,z).$
By using the relations (13), we get the second prolongation of the vector
field $Y$, which is defined by the next functions
$$
\Phi ^u  =  \phi _xx_u+\phi _yy_u+\phi _zz_u,\;\;
\Phi ^v  =  \phi _xx_v+\phi _yy_v+\phi _zz_v,
$$
$$
\begin{array}{ccl}
\Phi ^{uu} & = & \phi _{xx}x^2_u+\phi _{yy}y^2_u+\phi _{zz}z^2_u+
2\phi _{xy}x_uy_u+2\phi _{xz}x_uz_u+2\phi _{yz}y_uz_u+ \\
\noalign{\medskip }
&+&\phi _xx_{uu}+\phi _yy_{uu}+\phi _zz_{uu},  
\end{array}
$$
$$
\begin{array}{ccl}
\Phi ^{uv} & = & \phi _{xx}x_ux_v+\phi _{xy}x_vy_u+\phi _{xz}x_vz_u+
\phi _{xy}x_uy_v+\phi _{yy}y_uy_v+\phi _{yz}y_vz_u+ \\
\noalign{\medskip }
&+&\phi _{xz}z_vx_u+\phi _{yz}y_uz_v+\phi _{zz}z_uz_v+\phi _xx_{uv}+
\phi _yy_{uv}+\phi _zz_{uv}, 
\end{array}
$$
$$
\begin{array}{ccl}
\Phi ^{vv} & = & \phi _{xx}x^2_v+2\phi _{xy}x_vy_v+2\phi _{xz}x_vz_v+
2\phi _{yz}y_vz_v+\phi _{yy}y^2_v+\phi _{zz}z^2_v+ \\
\noalign{\medskip }
&+&\phi _xx_{vv}+\phi _yy_{vv}+\phi _zz_{vv},
\end{array} 
$$
and also the functions
$\Lambda ^u,\Lambda ^v,\Lambda ^{uu},\Lambda ^{uv},\Lambda ^{vv},
\Psi ^u,\Psi ^v,\Psi ^{uu},\Psi ^{uv},\Psi ^{vv}$ which are
analogously written by substituting $\phi $ with $\lambda $, and respectively
$\psi $. The PDEs system (15)+(16) is of maximal rank.
The infinitesimal invariance condition (14) for (15) turns in
$$
\left\{
\begin{array}{cc}
a\Phi ^u+b\Phi ^v-\Phi ^{uu} = 0 \\
h\Phi -\Phi ^{uv}=0 \\
a''\Phi ^u+b''\Phi ^v-\Phi ^{vv}=0 \\
.....................................\;.
\end{array} \right. \leqno(17)
$$
Let consider the first relation and let substitute the functions
$\Phi ^u,\Phi ^v$ and $\Phi ^{uu}$ gives by the above relations. We find
$$
a(\phi _xx_u+\phi _yy_u+\phi _zz_u)+b(\phi _xx_v+\phi _yy_v+
\phi _zz_v)-\phi _{xx}x^2_u-2\phi _{xy}x_uy_u-2\phi _{xz}x_uz_u-
$$
$$
-2\phi _{yz}y_uz_u-\phi _{yy}y^2_u-\phi _{zz}z^2_u-\phi _xx_{uu}-
\phi _yy_{uu}-\phi _zz_{uu}=0.
$$
We eliminate any dependencies among the derivatives of the $x,y,z$
by substituting
$$
x_{uu}  =  ax_u+bx_v,\;\;
y_{uu}  = ay_u+by_v,\;\;
z_{uu}  =  az_u+bz_v.
$$
Then it results
$$
\phi _{xx}x^2_u+\phi _{yy}y^2_u+\phi _{zz}z^2_u+2\phi _{xy}x_uy_u+
2\phi _{xz}x_uz_u+2\phi _{yz}y_uz_u=0.
$$
Now we equate the coefficients of the remaining unconstrained partial
derivatives of $x,y,z$ to zero, and we get the PDEs system
$$
\begin{array}{cccccc}
\phi _{xx}=0 & \phi _{yy}=0 & \phi _{zz}=0 & \phi _{xy}=0 & \phi _{yz}=0 &
\phi _{xz}=0.
\end{array}
$$
It follows the solution $\phi (x,y,z)=a_{11}x+a_{12}y+a_{13}z+k$,
with $a_{11},a_{12},a_{13},k\in {\bf R}.$ 
By substituting the function $\phi $ in the next two relations of the 
system (17), we find $k=0$ and thus
$\phi (x,y,z)=a_{11}x+a_{12}y+a_{13}z,\;a_{11},a_{12},a_{13}\in {\bf R}.$
Analogously, using the next six relations of the system (17) we get
$\lambda (x,y,z)=a_{21}x+a_{22}y+a_{23}z,\;a_{21},a_{22},a_{23}\in {\bf R},
\;\;\psi (x,y,z)=a_{31}x+a_{32}y+a_{33}z,\;a_{31},a_{12},a_{13}\in {\bf R}.$
Using the criterion of infinitesimal invariance (14), for (16),
we get
$$
\phi (y_uz_v-z_uy_v)+\lambda (x_vz_u-x_uz_v)+\psi (x_uy_v-x_vy_u)+
\Phi ^u(zy_v-yz_v)+\Phi ^v(yz_u-zy_u)+
$$
$$
+\Lambda ^u(xz_v-zx_v)+\Lambda ^v(zx_u-xz_u)+\Psi ^u(yx_v-xy_v)+
\Psi ^v(xy_u-yx_u)=0.
$$
If we substitute the functions $\Phi ^u,\Phi ^v,....$, then it results
$$
(x_uy_v-x_vy_u)(\psi -x\psi _x-y\psi _y+z\phi _x+z\lambda _y)+
(x_vz_u-x_uz_v)(\lambda -x\lambda _x-z\lambda _z+y\phi _x+
$$
$$
+y\psi _z)+(y_uz_v-y_vz_u)(\phi -y\phi _y-z\phi _z+x\lambda _y+x\psi _z)=0.
$$
We eliminate any dependencies among the derivatives of $x,y,z$ by using the
relation (16) itself, and we find $\phi _x+\lambda _y+\psi _z=0$ or
equivalent
$a_{33}+a_{11}+a_{22}=0$. Thus, the next functions
$$
\begin{array}{ccl}
\phi (x,y,z)& = & a_{11}x+a_{12}y+a_{13}z \\
\lambda (x,y,z)&  = & a_{21}x+a_{22}y+a_{23}z\\
\psi (x,y,z) & = & a_{31}x+a_{32}y-(a_{11}+a_{22})z
\end{array}
$$
define the infinitesimal generator $Y$ of the symmetry subgroup $G_1$:
$$
Y=a_{11}\left(x{\frac{\partial }{{\partial x}}}-z{\frac{\partial }{{\partial z}}%
} \right)+a_{22}\left(y{\frac{\partial }{{\partial y}}}-z {\frac{\partial }{{%
\partial z}}}\right)+a_{12}y{\frac{\partial }{{\partial x}}}+ a_{13}z{\frac{%
\partial }{{\partial x}}}+ 
$$
$$
+a_{21}x{\frac{\partial }{{\partial y}}}+a_{23}z{\frac{\partial }{{\partial y}}}+
a_{31}x{\frac{\partial }{{\partial z}}}+a_{32}y{\frac{\partial }{{\partial z}}}. 
$$

{\bf Theorem 4}.
{\it The Lie algebra associated to the subgroup $G_1$ of
the full symmetry group $G$ of the PDEs system (15)+(16) ($G_1$ acts on the space
of the dependent variables) is generated by the vector fields}
$$
Y_1=x{\frac{\partial }{{\partial x}}}-z{\frac{\partial }{{\partial z}}},\;\;
Y_2=y{\frac{\partial }{{\partial y}}}-z{\frac{\partial }{{\partial z}}},\;\;
Y_3=y{\frac{\partial }{{\partial x}}},\;\;Y_4=z{\frac{\partial }{{\partial x}%
}}\leqno (18)
$$
$$
Y_5=x{\frac{\partial }{{\partial y}}},\;\;
Y_6=z{\frac{\partial }{{\partial y}}},\;\;
Y_7=x{\frac{\partial }{{\partial z}}},\;\;Y_8=y{\frac{\partial }
{{\partial z}}},
$$
{\it and, thus the Lie subgroup $G_1$ is the unimodular subgroup of the
group of centroaffine transformations.}

Using this result, we can find group-invariant solutions
of the PDEs system (15)+(16).
For example, if we consider the subalgebra described by the vector fields
$Y_1$ and $Y_2$, then the function $F$ which is invariant under the associated
group, satisfies $Y_1(F)=0,$ and $Y_2(F)=0$. It results $F=\varphi (u,v,xyz)$
and the group-invariant solutions are defined by 
$u = C_1,\;v=C_2$ and $xyz=C_3$. Thus one gets the known \c{T}i\c{t}eica
surfaces
$$
z={C\over {xy}},\;\;C\in {\bf R}.\leqno(19)
$$
\par
\noindent
{\bf 3.2}.
We shall study if there is a subgroup $G_2$ of the symmetry group $G$
of the PDEs system (15)+(16) which acts on the space of the
independent variables $u,v$ of the system. Let suppose that it is and let $Z$ be its
infinitesimal generator. The vector field $Z$ is 
defined by the functions
$$
\zeta =\zeta (u,v),\;\;\eta =\eta (u,v),\;\;\phi =0,\;\;\lambda =0,
\;\;\psi =0.
$$
In this case, by using the above algorithm for finding the associated
symmetries, we get

{\bf Theorem 5}.
{\it The general vector field of the algebra of the infinitesimal symmetries
associated to the subgroup $G_2$, where $G_2$ is the subgroup of the full
symmetry group $G$ of the PDEs system (15), which acts on the space of the
independent variables, is
$$
Z=\zeta (u){\partial \over {\partial u}}+\eta (v){\partial \over {\partial v}},
$$
where the functions $\zeta $ and $\eta $ satisfy the next PDEs system:
$$
\left\{
\begin{array}{ccl}
& & \zeta a_u+\eta a_v+a\zeta _u+\zeta _{uu}=0 \\
& & \zeta b_u+\eta b_v-b\eta _v+2b\zeta _u=0 \\
& & \zeta h_u+\eta h_v+h(\zeta _u+\eta _v)=0 \\
& & \zeta a''_u+\eta a''_v-a''\zeta _u+2a''\eta _v=0 \\
& & \zeta b''_u+\eta b''_v+b''\eta _v+\eta _{vv}=0, 
\end{array} \right.
\leqno(20)
$$
and the functions $a,b,h,a'',b''$ satisfy the
integrability conditions (6)}.

One considers the cases:
\par
\noindent
1. If $\Sigma $ is a ruled \c{T}i\c{t}eica surface (7), then the
completely integrability conditions (6) are
$$
a={h_u\over {h}},\;\;
b={\varphi (u)\over {h}},\;\;
a''=0,\;\;b''={h_v\over {h}},
$$
where $h$ is a solution of the Liouville-\c{T}i\c{t}eica PDE (8).
In this case, the relations (20) turn in
$$
\left\{
\begin{array}{ccl}
\zeta h_u+\eta h_v+h(\zeta _u+\eta _v)=0 \\
\zeta ^3={k\over {\varphi }}, \\
hh_{uv}-h_uh_v=h^3.
\end{array} \right.
$$
Let us consider the change of variables
$\zeta ={1\over {U'}}$ and $\eta =-{1\over {V'}}$,
where $U=U(u)$ and $V=V(v)$. Then the first PDE implies
$h=U'V'\mu (U+V)$ and by substituting this function in the last PDE of
the above system (Liouville-\c{T}i\c{t}eica PDE), we find the following ODE
$$
\mu \mu ''-\mu '^2=\mu ^3,
$$
for which
$$
\mu (t)=
\left\{
\begin{array}{ccc}
{2\over {(t+C)^2}}, & k=0 \\
\\
l^2\over {2\hbox{cos} ^2({l\over {2}}t+C)}, & k=-l^2 \\
\\
{l^2 \over {2\hbox{sh}^2({l\over {2}}t+C)}}, & k=l^2,\;\;l>0.
\end{array}
\right.
$$
is the general solution, with $t=U+V$.
One substitutes in $h=U'V'\mu (U+V)$
and one consider the special change of the functions $\tilde U=F(U),\;
\tilde V=G(V)$, $\tilde U=U+C,\;\tilde V=V$: for $k=0$,
$\tilde U=\hbox{th} {l \over {2}}(U+C),\;\tilde V=\hbox{th}{l\over {2}}V$,
for $k=l^2$ and $\tilde U=\hbox{ctg}({l\over {2}}U+C),\;
\tilde V=\hbox{tg} {l\over {2}}V$, for $k=-l^2$, it results the general
solution of the Liouville-\c{T}i\c{t}eica ([15],[24]) namely,
$$
h(u,v)={2\tilde U'\tilde V'\over {(\tilde U+\tilde V)^2}}.
$$

\par
\noindent
2. If $\Sigma $ is a \c{T}i\c{t}eica surface which are not ruled
surface (9), then the
completely
integrability conditions (6) turn in
$$
a={h_u\over {h}},\;\;
b=a''={1\over {h}},\;\; 
b''={h_v\over {h}},
$$
where $h$ is a solution of the \c{T}i\c{t}eica PDE (10).
If we substitute these functions in the system (20), then we get
$$
\zeta _u=0,\;\;\;\;\eta _v=0,\;\;\;\;\zeta h_u+\eta h_v+h(\zeta _u+\eta _v)=0.
$$
It results $\zeta =C_1$, $\eta =C_2$ and $h=\mu (C_1v-C_2u)$ and thus
$$
Z=C_1{\partial \over {\partial u}}+C_2{\partial \over {\partial v}}.
$$
Let consider the \c{T}i\c{t}eica PDE (10) and let substitute the above function
$h=\mu (C_1v-C_2u)$.
We get the following ODE:
$$
-C_1C_2(\mu \mu ''-\mu '^2)=\mu ^3-1.
$$
\par
\noindent
a. If $C_1C_2=0$, then $\mu =1$ and $h=1$. This is the 
\c{T}i\c{t}eica solution [24].
\par
\noindent
b. If $C_1C_2\neq 0$, then denote $k=-{1\over {C_1C_2}}$. The above ODE turns
in
$$
\mu \mu ''-\mu '^2=k(\mu ^3-1).
$$
We can consider $k=1$. This ODE can be reduced to the following
$$
\mu '^{2}=2\mu ^3+C\mu ^2+1,\;\;C\in {\bf R},
$$
and using the change of function $\mu ={1\over {2}}g$, this becomes
$$
g'^2=g^3+Cg^2+4.
$$
Let $\lambda $ the real solution of the right side polynom of the above ODE.
It results that $\lambda \neq 0$, $\lambda $ is not a triple solution
of this and the ODE can writen as
$$
g'^2=(g-\lambda )\left(g^2-{4\over {\lambda ^2}}g-{4\over {\lambda }}\right).
$$
2.1. If $\lambda =-1$, then $C=-3$ and the ODE is
$$
g'^2=(g+1)(g-2)^2.
$$
If we consider $g={1\over {w^2}}+2$ than we have the ODE
$$
w'^2={1\over {4}}\left(3w^2+1\right),
$$
for which the solution $w=w(t),\;t=u+v$ is
$$
w(t)={1\over {\sqrt{3}}}\hbox{sh}\left({t\sqrt{3}\over {2}}+C_1\right),\;\;C_1\in
{\bf R}.
$$
We get the next solution of the \c{T}i\c{t}eica PDE: $h={1\over {2w^2}}+
1$ which in the case $C_1=0$, it is
$$
h(t)={3\over {2\hbox{sh}^2\left({t\sqrt{3}\over {2}}\right)}}+1,\;\;\;t=u+v.
\leqno(21)
$$
\par
\noindent
2.2. If $\lambda \neq -1$, then the right side polynom
has three distinct real solutions ($\lambda >-1$ or
$C<-3$) and repectively one is real and two complex ($\lambda <-1$
or $C>-3$). In this case, the integral
$$
J=\int {dg\over \sqrt{(g-\lambda )\left(g^2-{4\over {\lambda ^2}}g-{4\over
{\lambda }}\right)}}
$$
can be reduced to a
first genus elliptical integral [11]
$$
J=\int {d\varphi \over {\sqrt{1-k^2\sin \varphi }}}.
$$
We get: for $C\neq -3$, the solutions of the PDE
\c{T}i\c{t}eica with the form $h=\mu (u+v)$ are gives in the terms of
elliptical functions.

{\bf Proposition 1}.
{\it The solution (21) gives a revolution \c{T}i\c{t}eica surface.
Moreover, it is an associated ruled \c{T}i\c{t}eica surface}.

{\bf Proof}.
\^In [24], 164-174, \c{T}i\c{t}eica studied the revolution
surfaces defined by the system (11) and
make the remark that, in this case,
the function $h$ must satisfy
$h_u=h_v$, and thus $h=\mu (u+v)$.
\c{T}i\c{t}eica obtained the above ODE
$$
\mu \mu''-\mu '^2=\mu ^3-1,
$$
and not integrate it. He proved: by using
$$
{\mu '^2-2\mu ^3-1\over {4\mu ^2}}=-k^2,
$$
one finds the solution of the studied system: for $k\neq 0$:
$$
\theta (u,v)=k_1e^{\int {h'-1 \over {2h}}d\alpha }\cos k\beta+
k_2e^{\int {h'-1 \over {2h}}d\alpha }\sin k\beta +
k_3e^{\int {h^2 \over {h'+1}}d\alpha },
\leqno(22)
$$
and for $k=0$:
$$
\theta (u,v)=e^{\int {h'-1 \over {2h}}d\alpha }
\left[k_1\left(\beta ^2+\int {4\mu \over {\mu '+1}}d\alpha \right)+
k_2\beta +k_3\right],
$$
where $\alpha =u+v,\;\beta =u-v$ \c{s}i $k_1,k_2,k_3\in {\bf R}$.

On the other hand, ours calculus imply that 
$k^2=-{C\over {4}}$, and thus
the function (21) defines a revolution surface (22).
If we consider
$$
\tilde U=\hbox{th} {\sqrt{3}\over {2}}(U+C_1),\;\;
\tilde V=\hbox{th} {\sqrt{3}\over {2}}V,
$$
it results that the function $h$  can be written in the next form
$$
h(u,v)={2\tilde U'\tilde V'\over {(\tilde U+\tilde V)^2}}+1=
H(u,v)+1.
$$
But, the above calculus, implies that the function 
$H$ is a solution of Liouville-\c{T}i\c{t}eica PDE (8) and
this defines a ruled \c{T}i\c{t}eica surface.

{\bf Proposion 2}.
{\it The solution of \c{T}i\c{t}eica (22) is invariant
under the transformations subgroup of $G_2$, for which
the Lie algebra is generated (in the case of
not ruled \c{T}i\c{t}eica surface) by}
$$
Z=C_1{\partial \over {\partial u}}+C_2{\partial \over {\partial v}}.
$$
\par
\noindent
{\bf 3.3}.
In this part, we shall apply the symmetry group theory
for the PDEs system (11)
with the integrability conditions (6). Let $D\times \bar U^{(2)}$ be the
second order jet space associated to the PDEs system (11), whose
coordinates are the independent variables $u,v$, the dependent variable
$\theta $ and the derivatives of the dependent variable
till the order two.
Denote by $x^1=u,\;x^2=v$ and $u^1=\theta $ (in the second section).
Let consider $\bar M\subset D\times \bar U$ an open set and let
$$
\bar X=\zeta {\partial \over {\partial u}}+\eta {\partial \over {\partial v}}+
\alpha {\partial \over {\partial \theta}}
$$
be the infinitesimal generator of the symmetry group $\bar G$ of the PDEs
system (11), where $\zeta ,\;\eta $ and $\alpha $ are functions of $u,\;v$
and $\theta $. The relations (13) imply the first and the second
prolongations of the vector field $\bar X$
$$
pr^{(1)}\bar X=\bar X+\alpha ^u{\partial \over {\partial \theta _u}}+
\alpha ^v {\partial \over {\partial \theta _v}},
$$
$$
pr^{(2)}\bar X=pr^{(1)}\bar X+\alpha ^{uu}{\partial \over {\partial
\theta _{uu}}}+\alpha ^{uv}{\partial \over {\partial \theta _{uv}}}
+\alpha ^{vv} {\partial \over {\partial \theta _{vv}}},
$$
where
$$
\begin{array}{ccl}
\alpha ^u & = & D_u(\alpha -\zeta \theta _u-\eta \theta _v)+
\zeta \theta _{uu}+\eta \theta _{uv}, \\
\alpha ^v & = & D_v(\alpha -\zeta \theta _u-\eta \theta _v)+\zeta
\theta _{uv}+\eta \theta _{vv},
\end{array}
$$
$$
\begin{array}{ccl}
\alpha ^{uu} & = & D_{uu}(\alpha -\zeta \theta _u-\eta \theta _v)+
\zeta \theta _{uuu}+\eta \theta _{uuv}, \\
\alpha ^{uv} & = & D_{uv}(\alpha -\zeta \theta _u-\eta \theta _v)+
\zeta \theta _{uuv}+\eta \theta _{uvv}, \\
\alpha ^{vv} & = & D_{vv}(\alpha -\zeta \theta _u-\eta \theta _v)+
\zeta \theta _{uvv}+\eta \theta _{vvv}. 
\end{array}
$$

We shall study if there is a subgroup $\bar G_1$ of the symmetry group
$\bar G$, which acts on the space of the dependent variable $\theta $.
Let suppose that it is and let
$$
\bar Y=\alpha {\partial \over {\partial \theta }},\;\;\alpha
=\alpha (\theta ),
$$
be its infinitesimal generator. In this case, the above algorithm implies  

{\bf Theorem 6}.
{\it The Lie algebra of the infinitesimal symmetries associated to
the subgroup $\bar G_1$ of the full symmetry group $\bar G$ of the PDEs
system (11) ($\bar G_1$ acts on the space of the dependent variable $\theta $),
is generated by the vector field}
$$
\bar Y_1=\theta {\partial \over {\partial \theta }}. \leqno(23)
$$

Analogously, we study the subgroup $\bar G_2$ of the symmetry group $\bar G$,
which acts on the space of the independent variables $u,v$ of the PDEs system
(11). Let
$$
\bar  Z=\zeta {\partial \over {\partial u}}+\eta 
{\partial \over {\partial v}},\;\;\zeta =\zeta (u,v),\;\;\eta =\eta (u,v),
$$
be the infinitesimal generator of it. It results

{\bf Theorem 7}.
{\it The general vector field of the algebra of the infinitesimal symmetries
associated to the subgroup $\bar G_2$ of the symmetry group $\bar G$ of the
PDEs system (11) ($\bar G_2$ acts on the space of the independent variables
$u,v$), is
$$
\bar Z=\zeta (u) {\partial \over {\partial u}}+
\eta (v) {\partial \over {\partial v}}, \leqno(24)
$$
where $\zeta $ and $\eta $ satisfy the relations (20)}.

{\bf Remark}. The subgroup $G_2$ and $\bar G_2$
are the same actions on the space of the independent variables $u,v$.

\par
\noindent
{\bf 3.4}.
Now we shall study the symmetries of
the Liouville-\c{T}i\c{t}eica PDE (8) and \c{T}i\c{t}eica PDE (10).
Let consider the PDEs
$$
\omega _{uv}=e^{\omega }, \leqno(8')
$$
and respectively
$$
\omega _{uv}=e^{\omega }-e^{-2\omega }, \leqno(10')
$$
which are equivalent to Liouville-\c{T}i\c{t}eica PDE,
and respectively \c{T}i\c{t}eica PDE, where $\ln h=\omega $.
We remark that these PDEs belong of the next class of
second order PDE, of maximal rank,
$$
\omega _{uv}=H(\omega ), \leqno(25)
$$
which was studied by Sophus Lie himself.
Also Pucci, Saccomandi, Mansfield have considered such equations.
One proves [18] the next result

{\bf Theorem 8.}
{\it If $\zeta =\zeta (u,v,\omega),\;\eta =\eta (u,v,\omega)$ and
$\phi =\phi (u,v,\omega)$ are the solutions of the PDEs system
$$
\zeta _v=0\;\;\;\zeta _\omega =0\;\;\;\eta _u=0\;\;\;\eta _\omega =0\;\;\;
\phi _{\omega \omega }=0\;\;\;\phi _{u\omega }=0\;\;\;\phi _{v\omega }=0
\leqno(26)
$$
$$
\phi _{uv}+(\phi _{\omega }-\zeta _u-\eta _v-\phi )H-H'\phi =0,
$$
where $H=H(\omega)$, then
$$
X=\zeta {\partial \over {\partial u}}+\eta {\partial \over {\partial v}}
+\phi {\partial \over {\partial \omega}}
$$
is the infinitesimal generator of the symmetry group associated to
one PDE of the form (25).}

In the case of PDEs (8') and (10'), we get the results

{\bf Theorem 9.}
{\it The general vector field which describes the algebra
of infinitesimal symmetries associated to the 
Liouville-\c{T}i\c{t}eica PDE (8$'$) is the following
$$
W=f{\partial \over {\partial u}}+g{\partial \over {\partial v}}-
(f'+g'){\partial \over {\partial \omega }}, \leqno(27)
$$
where $f=f(u)$ and $g=g(v)$.}

{\bf Theorem 10.}
{\it The vector fields which generate the Lie algebra
of infinitesimal symmetries associated with the PDE 
\c{T}i\c{t}eica (10') are}
$$
U_1=u{\partial \over {\partial u}}-v{\partial \over {\partial v}},\;\;
U_2={\partial \over {\partial u}},\;\;U_3={\partial \over {\partial v}}.
\leqno(28)
$$

{\bf Remark}.
If $\omega =f(u,v)$ is o solution of the \c{T}i\c{t}eica PDE (8'), then
the following functions
$$
\omega ^{(1)}=f(e^{\varepsilon }u,e^{-\varepsilon }v),\;\;
\omega ^{(2)}=f(u-\varepsilon, v),\;\;
\omega ^{(3)}=f(u,v-\varepsilon ),
$$
where $\varepsilon $ is a real number, are also solutions of the equation.

The adjoint representation of the symmetry group of the 
\c{T}i\c{t}eica PDE (8'), is given by the next table

\vspace{0,3cm}
\centerline{$\scriptstyle{Table\;1}$}
\vspace{0,1cm}
\centerline{\begin{tabular}{|c|c|c|c|}
\hline
${Ad}$ & ${U_1}$ & ${U_2}$ &
${U_3}$ \\
\hline
${U_1}$ & ${U_1}$ & ${e^{\varepsilon }U_2}$ &
${e^{-\varepsilon }U_3}$ \\
\hline
${U_2}$ & ${U_1-\varepsilon U_2}$ & ${U_2}$ &
${U_3}$  \\
\hline
${U_3}$ & ${U_1+\varepsilon U_3}$ &
${U_2}$ & ${U_3}$  \\
\hline
\end{tabular}}

\par
\vspace{0,3cm}

The 1-dimensional subalgebras, described by $U_2,\;U_3,\;U_2-U_3$,
implies the finding the following group-invariant solutions:
\par
\noindent
1. For $U_2$, and respectively $U_3$, the solution is $\omega =0$.
Thus it results the \c{T}i\c{t}eica solution $h=1$ of the PDE (8).
\par
\noindent
2. In the case of the vector field $U_2-U_3$, the group-invariant solutions 
have the form $\omega =f(u+v)$ (respectively $h=\mu (u+v)$ for PDE (10)).
But this case was considered in the above section.

Now we look for the converse  of the Theorem 10:

{\bf Theorem 11}.
{\it The second order PDE invariant with respect to
the symmetry group (28) of \c{T}i\c{t}eica PDE,
has the form}
$$
H(\omega ,\omega _u\omega _v,\omega _{uv},\omega _{uu}\omega _{vv})=0.
\leqno(29)
$$

{\bf Proof.}
One considers the maximal chain of Lie subalgebras of the Lie algebra
associated with the studied group:
$$
\{U_2\}\subset \{U_2,U_3\}\subset \{U_1,U_2,U_3\}.
$$
and a second order PDE: $F(u,v,\omega ^{(2)})=0$ which is invariant
with resect this group. So it satisfies the criterion of infinitesimal
invariance for each second prolongation of these vectors fields.
\par
\noindent
1. For $U_2$, it results $pr^{(2)}U_2(F)=0,$ and thus
$F=F_1(v, \omega ^{(2)})$.
\par
\noindent
2. For $U_3$, the condition $pr^{(2)}U_3(F)=0,$
implies $F=F_2(\omega ^{(2)})$.
\par
\noindent
3. For $U_1$, the relation $pr^{(2)}U_1(F)=0$ is equivalent to
$$
U_1(F_2)-\omega _u{\partial F_2\over {\partial \omega _u}}+
\omega _v{\partial F_2\over {\partial \omega _v}}-2\omega _{uu}{\partial F_2
\over {\partial F_2\omega _{uu}}}+2\omega _{vv}{\partial F_2\over
{\partial \omega _{vv}}}=0,
$$
and finally, we get $F_2=H$ in (29).

\section{Lagrangians Associated to \c{T}i\c{t}eica PDEs}

{\bf 4.1}.
Our present study is for the inverse problem for the PDEs (8') and (10').
We recall that, the simple form of the inverse problem for a PDE
in the calculus of variations is to determine if this
is identically to an Euler-Lagrange PDE [1],[8],[16],[20],[21], [28],[31].

Let us consider the second order PDE
$$
\Delta (u,v,\omega ^{(2)})=0,\leqno(30)
$$
where $\omega ^{(2)}$ is the second prolongation of the unknown function
$\omega =\omega (u,v)$. The PDE (30) is {\it identically to
an Euler-Lagrange equation} if and only if the inegrability Helmholtz
conditions
$$
\left\{
\begin{array}{ccl}
{\frac{\partial \Delta }{{\partial \omega _u}}} & = & D_u\left(
{\frac{\partial \Delta }{{\partial \omega _{uu}}}}\right)+D_v\left({\frac{1}{{2}}}
{\frac{\partial \Delta }{{\partial \omega _{uv}}}} \right) \\ {\frac{\partial T}
{{\partial \omega _v}}} & = & D_u\left({\frac{1}{{2}}}{\frac{\partial \Delta }
{{\partial \omega _{uv}}}}\right)+D_v\left({\frac{\partial \Delta }{{\partial
\omega _{vv}}}} \right),
\end{array}
\right.\leqno(31) 
$$
are satisfied. In this case, there exists a function $L$,
called {\it Lagrangian}, such that the Euler-Lagrange PDE
$$
E(L)={\frac{\partial L}{{\partial \omega }}}-D_u\left({\frac{\partial L}{{%
\partial \omega _u}}} \right)-D_v\left({\frac{\partial L}{{\partial \omega _v}}}%
\right)=0 
$$
is equivalent to the PDE (30), in the sense
that every solution of the equation (30) is a solution of the
Euler-Lagrange equation $E(L)=0$ and conversely.

Also, the PDE (30) is called {\it equivalent to an Euler-Lagrange
equation} if there exists a nonzero function
$f=f(u,v,\omega ,\omega _u,\omega _v)$
such that $f\cdot \Delta =E(L)$. The function $f$ is called {\it variational
integrant factor}.

{\bf Theorem 12.}
{\it The Liouville-\c{T}i\c{t}eica PDE (8') and \c{T}i\c{t}eica PDE (10')
are Euler-Lagrange equations with the next associated Lagrangians 
$$
L_1(u,v,\omega ^{(1)})=-{1\over {2}}\omega _u\omega _v-e^{\omega },
\leqno(32)
$$
and}
$$
L_2(u,v,\omega ^{(1)})=-{1\over {2}}\omega _u\omega _v-e^{\omega }
-{1\over {2}}e^{-2\omega }. \leqno(33)
$$

{\bf Proof}.
Indeed, one verifies that the Helmholtz integrability conditions (31)
are satisfied and one verifies that $L_1$ and $L_2$ can be considered
associated Lagrangians.

{\bf Remark}.
The PDEs (8) and (10) are equivalent to Euler-Lagrange PDEs, with
the variational integrant factor ${1\over {h^3}}$.

\par
\noindent
{\bf 4.2.}.
We make a short presentation of the theory of variational symmetry groups
for the functionals 
$$
{\cal L}[\omega ]=\int _{}^{}\int _{\Omega _0}^{}L(u,v,\omega ^{(1)})dudv,
\leqno(34)
$$
with $\Omega _0$ is a domain in {\bf R$^2$} [20],[21].

Let $D\subset \Omega _0$ be a subdomain, $U$ an open set in {\bf R}
and $M\subset D\times U$ an open set. We consider
$\omega \in C^2(D),\;\omega =f(u,v)$
such that $\Gamma _\omega =\{(u,v,\omega (u,v))\vert (u,v)\in D\}\subset M$.
A local group $G$ of transformations on $M$ is called {\it
variational symmetry group for the functional} (34),
if $g_{\varepsilon }\in G,\;g_{\varepsilon }(u,v,\omega )= (\bar
u,\bar v,\bar \omega )$, then the function $\bar \omega =\bar f(\bar u,\bar v)=
(g\cdot f)(\bar u,\bar v)$ is defined on $\bar \Omega \subset \Omega _0$ and
$$
\int _{}^{}\int _{\bar D}^{}L(\bar u,\bar v,pr^{(1)} \bar f(\bar u,\bar
v))d\bar u d\bar v= \int \int _{D}L(u,v,pr^{(1)}f(u,v))dudv. 
$$

{\bf Theorem 13 (Infinitesimal criterion for the variational problem)}.
$\;\;$
{\it A connected group $G$ of transformations acting on $M\subset \Omega_0
\times U$ is a group of variational symmetries for the functional (34) if
and only if
$$
\hbox{pr}^{(1)}X(L)+L\;Div\xi =0, \leqno(35) 
$$
is satisfied for $\forall (u,v,\omega ^{(2)})\in M^{(2)}\subset
D\times U^{(2)}$ and for any infinitesimal generator  
$$
X=\zeta (u,v,\omega ){\frac{\partial }{{\partial u}}}+\eta (u,v,\omega ) {\frac{\partial 
}{{\partial v}}}+\phi (u,v,\omega ){\frac{\partial }{{\partial \omega }}} 
$$
of $G$, where $\xi =(\zeta ,\eta )$ and $Div\xi =D_u\zeta +D_v\eta $}.

{\bf Theorem 14.}
{\it If $G$ is a variational symmetry group of the functional (34),
then $G$ is a symmetry group of Euler-Lagrange equation $E(L)=0$}.

The converse of Theorem 14 is generally false.

Let us consider the PDE (30).
{\it A conservation law} is a divergence expression $Div\;P=0$
wich vanishes for all solutions $u=f(x)$ of the given PDE.
Here $P=(P^1,P^2)$ with $Div\;P=D_uP^1+D_vP^2,$ {\it the total divergence.}
The function $P^1$ is called {\it
flow associated}
and $P^2$ is called {\it conserved density to the conservation law}.
It results that
there exists a function $Q$ such that 
$$
Div\;P=Q\cdot \Delta .\leqno(36) 
$$
This relation is called {\it the characteristic form of the conservation
law}, and $Q$ is called {\it the characteristic of the conservation law.}

Let
$$
X=\zeta (u,v,\omega ){\frac{\partial }{{\partial u}}}+\eta (u,v,\omega ){\frac{\partial 
}{{\partial v}}}+\phi (u,v,\omega ){\frac{\partial }{{\partial \omega }}} 
$$
be a vector field on $M$. The vector field  
$$
X_Q=Q{\frac{\partial }{{\partial u}}},\;\;Q=\phi -\zeta \omega _u-\eta \omega _v, 
$$
is called {\it vector field of evolution associated to $X$},
and $Q$ is called {\it the characteristic associated to $X$}.

{\bf Theorem 15 (Noether Theorem)}.
{\it Let $G$ be a local Lie group of transformations, which is a symmetry
group of the variational problem (34) and let
$$
X=\zeta (u,v,\omega ){\frac{\partial }{{\partial u}}}+\eta (u,v,\omega ){\frac{\partial 
}{{\partial v}}}+\phi (u,v,\omega ){\frac{\partial }{{\partial \omega }}} 
$$
be the infinitesimal generator of $G$. The characteristic $Q$ of the field $X$
is also a characteristic of the conservation law for the associated 
Euler-Lagrange equation $E(L)=0$.}

One proves ([20], 356) that for the Lagrangian
$L=L(u,v,\omega ^{(1)})$ we have
$$
P=-(A+L\xi )=-(A^1+L\zeta ,A^2+L\eta )=(P^1,P^2),\;\;A=(A^1,A^2),\leqno(37) 
$$
where  
$A^1=Q\cdot E^{(u)}(L),\;\;\;A^2=Q\cdot E^{(v)}(L)$. The operators 
$E^{(u)}(L)={\frac{\partial L}{{\partial \omega _u}}}$ and
$E^{(v)}(L)={\frac{\partial L}{{\partial \omega _v}}}$
are called {\it first order Euler operators.}

\par
\noindent
{\bf 4.3}.
Let us consider the first order Lagrangians (32) and (33) and the
associated functionals
$$
{\cal L}[\omega ]=\int \int _{D}^{}
L_1(u,v,\omega ^{(1)})dudv \leqno(38)
$$
and
$$
\bar {\cal L}[\omega ]=\int \int _{D}^{}
L_2(u,v,\omega ^{(1)})dudv, \leqno(39)
$$
with $D$ is a domain in {\bf R}$^2$ and $\omega \in C^{2}(D)$.

{\bf Theorem 16}.
{\it The Lie algebras of the variational symmetry groups of the functional
(38) is described by the vector fields}
$$
W_1=u{\partial\over {\partial u}}-{\partial \over {\partial \omega }},
\;\;\;
W_2=v{\partial\over {\partial v}}-{\partial \over {\partial \omega }},
\;\;\;
W_3={\partial\over {\partial u}},\;\;\;W_4={\partial \over {\partial v}},
\leqno(40)
$$

{\bf Proof.} According with Theorem 14, the vector fields which determine
the Lie algebra of the variational symmetry group are founded
between the vector fields of the Lie algebra of the symmetry group
of the associated Euler-Lagrange equation.
The condition (35) must be verified only for the vector fields in
the algebra of the symmetry group of PDE (8$'$).
Let us consider the vector field 
$$
W=f{\partial \over {\partial u}}+g{\partial \over {\partial v}}-
(f'+g'){\partial \over {\partial \omega }},\;\;f=f(u),\;g=g(v),
$$
which is given by the relation (27) and let 
$$
pr^{(2)}W=W-\left(f''+f'\omega _u\right){\partial
\over {\partial \omega _u}}-\left(g''+g'\omega _v\right){\partial
\over {\partial \omega _v}},
$$
be the second prolongation. Introducing $\xi =(f,g)$ and $Div\xi =f'+g'$
in the relation (35), this turns in
$f''\omega _v+g''\omega _u=0$ and it implies $f''=g''=0$. Thus, we get
$f  = C_1u+C_3,\;g  =  C_2v+C_4$ and also
$$
W=(C_1u+C_3){\partial \over {\partial u}}+
(C_2v+C_4){\partial \over {\partial v}}-(C_1+C_2)
{\partial \over {\partial \omega }}=
$$
$$
=C_1\left(u{\partial \over {\partial u}}-{\partial \over {\partial \omega }}
\right)+C_2\left(v{\partial \over {\partial v}}-{\partial \over {\partial
\omega }}\right)+C_3{\partial \over {\partial u}}+C_4{\partial \over
{\partial v}}.
$$

{\bf Theorem 17}.
{\it The Lie algebras of the variational symmetry groups of the functional
(39) is described by the vector fields}
$$
U_1=u{\partial\over {\partial u}}-v{\partial\over {\partial v}},\;\;\;
U_2={\partial\over {\partial u}}\;\;\;U_3={\partial \over {\partial v}}.
\leqno(41)
$$

{\bf Proposition 3}.
{\it The
associated flows and respectively the
conserved density in the case of Liouville-\c{T}i\c{t}eica PDE (8')
and respectively
for \c{T}i\c{t}eica PDE (10') are}

\vspace {0,3cm}
\centerline{$\scriptstyle{Table\;2}$}
\vspace {0,1cm}
\centerline{\begin{tabular}{|c|c|c|}
\hline
${-W_i}$  & ${P^1}$ & ${P^2}$ \\
\hline
${-W_1}$ & ${{1\over {2}}\omega _v-ue^{\omega }}$ &
${{1\over {2}}\omega _u(1+u\omega _u)}$ \\
\hline
${-W_2}$ & ${{1\over {2}}\omega _v(1+v\omega _v)}$ &
${{1\over {2}}\omega _u-ve^{\omega }}$ \\
\hline
${-W_3}$ & ${-e^{\omega }}$ &
${{1\over {2}}\omega ^2_u}$ \\
\hline
${-W_4}$ & ${{1\over {2}}\omega ^2_v}$ &
${-e^{\omega }}$ \\
\hline
\end{tabular}}
\par

\vspace{0,5cm}

\centerline{$\scriptstyle{Table\;3}$}
\vspace{0,1cm}

\centerline{\begin {tabular}{|c|c|c|}
\hline
${-U_i}$  & ${P^1}$ & ${P^2}$ \\
\hline
${-U_1}$ & ${-{1\over {2}}ue^{-2\omega }-
{1\over {2}}v\omega ^2_v-ue^{\omega }}$ &
${{1\over {2}}u\omega ^2_u+ve^\omega +{1\over {2}}ve^{-2\omega }}$ \\
\hline
${-U_2}$ & ${-e^\omega -{1\over {2}}e^{-2\omega }}$ &
${{1\over {2}}\omega ^2_u}$ \\
\hline
${-U_3}$ & ${{1\over {2}}\omega ^2_v}$ &
${-e^\omega -{1\over {2}}e^{-2\omega }}$ \\
\hline
\end{tabular}}
\par
\vspace{0,3cm}

{\bf Proof}.
For example, the caracteristic associated to the vector field $-W_3$ is
$Q=\omega _u$. Replacing in the relations (37), we obtain
$A^1=-{1\over {2}}\omega _u\omega _v,\;\;A^2=-{1\over {2}}\omega ^2_u$
and thus $P^1=-e^{\omega },\;\;P^2={1\over {2}}\omega ^2_u$.

\par
\centerline{University "Politehnica" of Bucharest}
\centerline{Departament of Mathematics I}
\centerline{Splaiul Independen\c{t}ei 313}
\centerline{77206 Bucharest Romania}
\centerline{e-mail:nbila@mathem.pub.ro}
\end{document}